\pdfoutput=1    % arXiv recommends this as the first line; see http://arxiv.org/help/submit_tex

\documentclass[reqno, 12pt]{amsart}

\usepackage[T1]{fontenc}
\usepackage[latin1]{inputenc}
\usepackage{amsmath}
\usepackage{amssymb,verbatim}
\usepackage{amsfonts}
\usepackage{amsthm}
\usepackage{pifont}
\usepackage{wrapfig} % Paket zur Positionierung einbinden
\usepackage{graphicx}
\usepackage{enumitem}
\usepackage{tikz}

\usepackage{cases}
\usepackage{amsmath,mathtools}  % 5/18/2014; for dcases
\usepackage{hyperref}  % displays a bookmarks window, which lists section headings and References

\hypersetup{colorlinks=true, linkcolor=red, citecolor=red, pdfstartview=FitH, linktocpage=false}

\usepackage{comment}  % so we can comment out a block of text

\usepackage{parskip}  % skip a line between paragraphs

\newcommand{\beq}{\begin{equation*}}
\newcommand{\eeq}{\end{equation*}}
\newcommand{\beqn}{\begin{equation}}
\newcommand{\eeqn}{\end{equation}}
\newcommand{\ba}{\begin{array}}
\newcommand{\ea}{\end{array}}
\def\bal#1\eal{\begin{align*}\begin{split}#1\end{split}\end{align*}}
\def\baln#1\ealn{\begin{align}\begin{split}#1\end{split}\end{align}}
\newcommand{\tn}{\textnormal}
\newcommand{\dd}{\mathrm{d}}
\newcommand{\ee}{\mathrm{e}}
\newcommand{\si}{\mathrm{sinc}}

\usepackage[ top = 2.50cm, bottom = 2.50cm, left = 2.0cm, right = 2.0cm ]{geometry}

\usepackage{listings}
\usepackage{textcomp}  % required
\makeatletter
\def\lst@outputspace{{\ifx\lst@bkgcolor\empty\color{white}\else\lst@bkgcolor\fi\lst@visiblespace}}
\makeatother
\begin{document}
\title{Sinc integrals and tiny numbers}
\author{Uwe B\"asel\quad and\quad Robert Baillie}
%\date{\today,\;\currenttime\vspace{-0.5cm}}
%\date{\today}

\begin{abstract}
%A result of David and Jon Borwein is applied to evaluate a sequence of highly-oscillatory integrals whose integrands are the products of a rapidly growing number of sinc functions.
%The value of each integral is given in the form $\pi(1-t)/2$, where the numbers $t$ quickly become very tiny.
%Using the Euler-Maclaurin summation formula and Stirling's approximation, we calculate these numbers to high precision.
%For example, the integrand of the tenth integral in the sequence is the product of 68100152 sinc functions, and the value of the corresponding $t$ is approximately
We apply a result of David and Jon Borwein to evaluate a sequence of highly-oscillatory integrals whose integrands are the products of a rapidly growing number of sinc functions.
The value of each integral is given in the form $\pi(1-t)/2$, where the numbers $t$ quickly become very tiny.
Using the Euler-Maclaurin summation formula, we calculate these numbers to high precision.
For example, the integrand of the tenth integral in the sequence is the product of 68100152 sinc functions.
The corresponding $t$ is approximately
\[
9.6492736004286844634795531209398105309232 \cdot 10^{-554381308} \, .
\]
\end{abstract}

\maketitle

\section{Introduction} \label{S:Introduction}

Leonhard Euler knew that
\beqn \label{sinus}
  \int_0^\infty\frac{\sin x}{x}\:\dd x = \frac{\pi}{2}
\eeqn
at the latest by 1781 \cite[p.\ 324]{Remmert}, and there exist several proofs for it (see, for example \cite{Hardy}, \cite[p.~324]{Remmert}, and, for a proof due to Lobachevsky, \cite[pp.\ 635-636] {Fichtenholz}). The substitution $x\mapsto a_0 x$ with a real number $a_0$ immediately shows, more generally, that (see also \cite[pp.\ 83-84]{Nahin})
\beqn \label{sinus-a_0}
  I_0 = \int_0^\infty\frac{\sin(a_0 x)}{x}\:\dd x = \frac{\pi}{2}
  \quad\mbox{if}\quad a_0 > 0\,.
\eeqn
P\'olya \cite[pp.\ 208-209]{Polya} and D. \& J. Borwein \cite[pp.\ 78-79]{BorBor} derived in different ways the general evaluation of the integral  
\beqn \label{E:sin-prod}
  I_n = \int_0^\infty\frac{\sin(a_0 x)}{x}\,\frac{\sin(a_1 x)}{x}\cdots\frac{\sin(a_n x)}{x}\,\dd x
\eeqn
where $a_0,a_1,\ldots,a_n$ are real numbers. If, in addition, $a_0,a_1,\ldots,a_n > 0$ with
\beqn \label{E:condition}
  a_0 \ge s(n) := \sum_{k=1}^n a_k \,,
\eeqn
then the general solution simplifies to
\beqn \label{E:sin-prod-solution-1}
  I_n = \frac{\pi}{2}\,a_1\,a_2\cdots a_n \, ,
\eeqn
see \cite[p.\ 654]{Fichtenholz} and \cite[pp.\ 78-79]{BorBor}.
Furthermore, Corollary 1 of \cite{BorBor} (see also Theorem 2 of \cite{BBB}) says that, if
\beqn \label{E:2a_k)}
  2a_k \ge a_n > 0 \;\;\tn{for}\;\; k = 0,1,\ldots,n-1,
\eeqn
and $n$ is such that the sum of $a_1 + a_2 + \ldots + a_n$ first exceeds $a_0$,
\beqn \label{E:SumExceedsa0}
  s(n) > a_0 \ge s(n-1) \, ,
\eeqn
%then the general solution simplifies to
then we have this formula for the exact value of the integral:
\beqn \label{E:sin-prod-solution-2}
  I_n = \frac{\pi}{2} \left\{ \prod_{k=1}^n a_k - \frac{(a_1 + a_2 + \cdots + a_n - a_0)^n}{2^{n-1} \, n!} \right\} .
\eeqn

Now we consider the integral 
\beqn \label{E:J_n}
  J_n
= a_0\int_0^\infty\prod_{k=0}^n\si(a_k x)\,\dd x
= \frac{I_n}{a_1\,a_2\cdots a_n}
\eeqn
where the $\si$ function is defined as
\begin{equation*}
  \si(x) = 
  \begin{cases}
	\: \sin(x)/x & \text{if \, $x \ne 0$}\,,\\[0.1cm]
	\: 1 & \text{if \, $x = 0$} \, .
  \end{cases}
\end{equation*}
Eq.\ \eqref{sinus-a_0} yields
\beqn \label{E:J_0=pi/2}
  J_0 = a_0\int_0^\infty\si(a_0 x)\,\dd x = \frac{\pi}{2}
  \quad\mbox{if}\quad
  a_0 > 0\,.
\eeqn
For $n=1,2,\ldots$, and $a_0,a_1,\ldots,a_n>0$, from \eqref{E:condition}, \eqref{E:sin-prod-solution-1}, \eqref{E:2a_k)}, \eqref{E:SumExceedsa0}, \eqref{E:sin-prod-solution-2}, and \eqref{E:J_n} it follows that
\begin{numcases} 
 {J_n = }
 %\pi/2 \hspace{-0.3cm} & if $a_0 \ge s(n)$\,, \label{E:J_n=pi/2-1}\\[0.1cm]
 \frac{\pi}{2} \hspace{-0.3cm} & if $a_0 \ge s(n)$\,, \label{E:J_n=pi/2-1}\\[0.1cm]
 %\pi(1-t_n)/2 \hspace{-0.3cm} & if \eqref{E:2a_k)} and \eqref{E:SumExceedsa0} hold\,, \label{E:J_n-t_n}
 \frac{\pi}{2}(1-t_n) \hspace{-0.3cm} & if \eqref{E:2a_k)} and \eqref{E:SumExceedsa0} hold\,, \label{E:J_n-t_n}
\end{numcases}
where
\beqn \label{E:t_n}
 t_n = \frac{(a_1 + a_2 + \cdots + a_n - a_0)^n}{2^{n-1}\,n!\,\prod_{k=1}^n a_k}\,.
\eeqn

%\begin{equation*}
%\frac{\pi - x}{2} =
%  \begin{dcases}  % better vertical spacing than \begin{cases}
%    \sum_{n=1}^{\infty} \frac{\sin(nx)}{n}  &\text{for $0 < x < 2 \pi$,}\\
%    \sum_{n=1}^{\infty} \frac{\sin(nx)}{n} \cdot \frac{\sin(n)}{n}  &\text{for $1 \leq x \leq 2\pi - 1$.}\\
%  \end{dcases}
%\end{equation*}

Using Theorem 1 (ii) of \cite{BorBor}, we see that
\beqn \label{E:less_pi(2}
  0 < J_{n+1} \le J_n < \pi/2 \quad\tn{if}\quad
  a_{n+1} \le a_0 < s(n)\,,\;
  n \ge 1\,. 
\eeqn 
Schmid \cite[pp.\ 13-16]{Schmid} proves that $J_{n+1} < J_n <\pi/2$ if $\{a_k\}$ is a montonically non-increasing series of positive real numbers with $a_0<s(n)$. 

In our applications below, the $a_k$ are defined as
\beqn \label{E:a_k}
  \tn{$a_0$ is an integer $\ge 1$}, \qquad a_k = \frac{1}{2k-1}\;\;\tn{for}\;\; k=1,2,\ldots,n,
\eeqn
and we write $J_n(a_0)$ and $t_n(a_0)$ instead of $J_n$ and $t_n$, respectively. 

From \eqref{E:J_n=pi/2-1} and \eqref{E:less_pi(2} we know (see also \cite[pp.\ 3-4]{Baesel}, with a slight change in notation) that
\begin{align}
  J_n(a_0) = {} & \frac{\pi}{2} \quad \mbox{if} \quad a_0 \ge s(n) = \sum_{k=1}^n \frac{1}{2 k - 1}\,,\label{E:EqualsPiOver2}\\
  J_n(a_0) < {} & \frac{\pi}{2} \quad \mbox{if} \quad a_0 < s(n)\,.\label{E:LessThanPiOver2} 
\end{align}
It is easy to see that our $a_k$ as defined in \eqref{E:a_k} satisfy the inequalities \eqref{E:2a_k)}:
our $a_k$ are all positive with $a_k \ge a_n$ for $k = 0, 1, \ldots, n - 1$, which implies that $2 a_k \ge a_n$.
If, in addition, the inequalities \eqref{E:SumExceedsa0} are satisfied, then we are able to compute the \emph{exact} value of $J_n(a_0)$, given by Equations \eqref{E:J_n-t_n} and \eqref{E:t_n}.
$t_n(a_0)$ can be a \emph{very} tiny number.
{\em Our main aim in this paper is to show how to calculate these tiny numbers with high precision.}
%For given value of $a_0$, the first step in such a calculation is to calculate exactly the $n$ that satisfies \eqref{E:SumExceedsa0}.

% we did not use double factorials to evaluate the integrals for a0 up to 10.
% so, move the double factorials to the place where we show how to use euler-maclaurin to calculate the large factorials.

%Observing that 
%\beq
%  \prod_{k=1}^n a_k
%= \prod_{k=1}^n \frac{1}{2k-1}
%= \frac{1}{2n-1}\cdot\frac{1}{2n-3}\cdots\frac{1}{3}\cdot\frac{1}{1}
%= \frac{1}{(2n-1)!!}\,,
%\eeq
%we can rewrite Eq.\ \eqref{E:J_n-t_n} as
%\beqn \label{E:J_n-t_n-2}
%  J_n(a_0) = \frac{\pi}{2} \left\{ 1 - \frac{(s(n) - a_0)^n \, (2n-1)!!}{2^{n-1} \, n!} \right\} .
%\eeqn

With \textit{Mathematica} we can calculate a few of these integrals directly.
For example,
\beq
  J_2(1) = \int_0^\infty \si(x) \, \si \left( \frac{x}{1} \right) \, \si \left( \frac{x}{3} \right) \dd x
  = \frac{11}{24} \, \pi \approx 0.458333 \, \pi \, .
\eeq
In this example, we have $a_0 = 1$, $a_1 = 1$, and $a_2 = 1/3$.
Note that
\[
a_1 \leq a_0 < s(2) = a_1 + a_2 = \frac{4}{3} \, ,
\]
so using the inequalities \eqref{E:SumExceedsa0}, we have $n = 2$. Eq.\ \eqref{E:J_0=pi/2} tells us that
\[
  J_0(1)
  = \int_0^\infty \si(x) \, \dd x
  = \frac{\pi}{2}\,,
\]
and Eq.\ \eqref{E:J_n=pi/2-1} delivers
\[
  J_1(1)
  = \int_0^\infty \si(x) \cdot \si(x) \, \dd x
  = \frac{\pi}{2} \, .
\]
For the equation $J_1(1) = J_0(1)$ see, for example, \cite{Baillie6241}, \cite{BBB}, \cite{Boas}, \cite[p.\ 324]{Remmert}.
From Eq.\ \eqref{E:J_n-t_n} we get the already known result
\[
  J_2(1)
= \frac{\pi}{2} \left\{
% 1 - \frac{ (\frac{4}{3} - 1)^2\,3!! }{ 2^1 \cdot 2! }
 1 - \frac{ (\frac{4}{3} - 1)^2 }{ 2^1 \cdot 2! \cdot \frac{1}{1} \cdot \frac{1}{3} }  % removed !!
   \right\}
   = \frac{11 \pi}{24}\,.
\]

As we include more $\si$ functions in the integrand, it generally takes more time for \textit{Mathematica} to evaluate the integral.
\textit{Mathematica} is able to calculate that
\[
  J_n(2) = 2 \int_0^\infty \si(2x) \, \si \left( \frac{x}{1} \right) \cdots\, \si \left( \frac{x}{2n - 1} \right)\dd x
  = \frac{\pi}{2}
\]
for all $n = 1, 2, \ldots, 7$.

\textit{Mathematica} is also able to calculate (see also \cite[p.\ 4]{Baesel})
\begin{align*}
  J_8(2) = {} & 2 \int_0^\infty \si(2x) \, \si \left( \frac{x}{1} \right) \, \si \left( \frac{x}{3} \right) \cdots\, \si \left( \frac{x}{15} \right) \dd x \\[0.1cm]
  = {} & \frac{168579263752211300739165075916829279}
  {337158527504429357358419617830000000}\,\pi
  \approx 0.49999999999998998115\,\pi\,.
\end{align*}
%We can also evaluate this integral using Eq.\ \eqref{E:J_n-t_n-2}.
We can also evaluate this integral using Equations \eqref{E:J_n-t_n} and \eqref{E:t_n}.
For $J_8(2)$, we have $a_0 = 2$, $a_1 = 1$, $a_2 = 1/3$, \dots, $a_8 = 1/15$.
Then
\[
  \sum_{k = 1}^7 a_k = \frac{88069}{45045} = 1.95513\ldots < a_0 < 
  2.02181\ldots = \frac{91072}{45045} = \sum_{k = 1}^8 a_k \, ,
\]
and therefore,
\begin{align*}
  J_8(2) = {} & \frac{\pi}{2} \left\{
% 1 - \frac{(\frac{91072}{45045} - 2)^8\,15!!}{2^7 \cdot 8!} \right\}  % removed !!
 1 - \frac{ (\frac{91072}{45045} - 2)^8 }{ 2^7 \cdot 8! \cdot \frac{1}{1} \cdot \frac{1}{3} \cdots \frac{1}{15} } \right\}
  = \frac{\pi}{2} \left\{1 - \frac{3377940044732998170721}{168579263752214678679209808915000000} \right\} \\[0.2cm]
  = {} & \frac{168579263752211300739165075916829279}{337158527504429357358419617830000000} \, \pi \, .
\end{align*}

\section{First Calculations}

Given the value of $a_0$, the first task is to find the value of $n$ such that the inequalities in \eqref{E:SumExceedsa0} are satisfied.
For $a_0 \leq 10$, we can find the corresponding $n$ by brute force, that is, by simply computing partial sums of the $a_k$ until the sum exceeds $a_0$.
Then, with the help of \eqref{E:t_n} we compute the decimal approximations of the $t$ values for $a_0=1,2,\ldots,10$, shown in Table \ref{Ta:SmallTable}.
%The sums for Table \ref{Ta:SmallTable} were computed to 60 digits.
%The values of the integrals are rounded in the last ($40^{\text{th}}$) decimal place.
The $t$ values are rounded in the last ($40^{\text{th}}$) decimal place.
In \textit{Mathematica} on a standard laptop, only the last two $n$ values took more than a minute to calculate.
For $a_0 \geq 9$, we can calculate the $n$ and $t$ values much more quickly using the Euler-Maclaurin summation formula, as discussed in later sections of this paper.
One sees that the numbers $t_n(a_0)$ quickly become rather tiny.

\begin{table}[ht]
\begin{center}
\renewcommand{\arraystretch}{1.2}
\begin{tabular}{|r|r|l|} \hline
\rule{0pt}{14pt}
$a_0$ & \multicolumn{1}{c|}{$n$} & \multicolumn{1}{c|}{$t_n(a_0)$} \\[3pt] \hline\hline
  1 & $1$         &  $8.3333333333333333333333333333333333333333 \cdot 10^{-2}$ \\
  2 & $8$         &  $2.0037696034181553438737278689296078869394 \cdot 10^{-14}$ \\
  3 & $57$        &  $4.2814541036242680424608725308114449824436 \cdot 10^{-143}$ \\
  4 & $419$       &  $2.3710999975681168329914604542463318249729 \cdot 10^{-1326}$ \\
  5 & $3092$      &  $2.5899544469237193354708111256703110686749 \cdot 10^{-13544}$ \\
  6 & $22846$     &  $1.3063312175270580087816919230036297407401 \cdot 10^{-107025}$ \\
  7 & $168804$    &  $1.2084753305711806308265054034357601336007 \cdot 10^{-970071}$ \\
  8 & $1247298$   &  $6.9312222993226491135738066834549327656340 \cdot 10^{-8742945}$ \\
  9 & $9216354$   &  $3.3216970999036058275367686941671199966288 \cdot 10^{-67342884}$ \\
 10 & $68100151$  &  $9.6492736004286844634795531209398105309232 \cdot 10^{-554381308}$ \\ \hline
\end{tabular}
\vspace{0.2cm}
\renewcommand{\arraystretch}{1}
\caption{Values of $t_n(a_0)$ for the evaluation of the integrals $J_n(a_0) = \frac{\pi}{2} \, \{1-t_n(a_0)\}$}
\label{Ta:SmallTable}
\end{center}
\end{table}

Look what happens when we take ratios of successive $n$ values:
\begin{align*}
419/57           &\approx 7.35087719 \\
3092/419         &\approx 7.37947494 \\
22846/3092       &\approx 7.38874515 \\
168804/22846     &\approx 7.38877703 \\
1247298/168804   &\approx 7.38903107 \\
9216354/1247298  &\approx 7.38905538 \\
68100151/9216354 &\approx 7.38905548
\end{align*}
These ratios appear to be approaching $\ee^2 \approx 7.38905610$. That is, the $n$ that corresponds to $a_0+1$ is roughly $\ee^2$ times the $n$ that corresponds to $a_0$. Here is the explanation. The sum of $N$ terms of the harmonic series,
\[
  \sum_{k=1}^{N} \frac{1}{k} \, ,
\]
is about $\ln(N)$. We have $\ln(\ee \cdot N) = \ln(N) + 1$. Therefore, if $N$ terms of the harmonic series are required to reach a sum $S (\approx \ln(N))$, then about $\ee \cdot N$ terms are needed to make the sum reach $S + 1$.
The terms in our series,
\[
  \sum_{k=1}^{n} \frac{1}{2 k - 1} \, ,
\]
are about $1/2$ as large as the corresponding terms in the harmonic series. Therefore, to increase our sum by 1 requires about as many terms as the harmonic series needs to increase its sum by 2, which is about $\ee\cdot\ee = \ee^2$.

As noted above, Eq. \eqref{E:J_n-t_n} gives the {\em exact} value of the integral. The expression on the right side of Eq.\ \eqref{E:J_n-t_n} may be written as
\[
\frac{\pi}{2}\left(1-\frac{P}{Q}\right),
\]
where $P$ and $Q$ are integers. As $a_0$ increases, $P$ and $Q$ quickly become very large. For example, with $a_0=6$, $P$ and $Q$ have 453130185 and 453237210 digits, respectively. Displaying the first and last 20 digits for this case, we have
\[
\frac{P}{Q} =
  \frac{34293043773392420460 \text{ (453130145 digits) } 34573721229967337961}
       {26251415654224851611 \text{ (453237170 digits) }  00000000000000000000} \,.
\]

\section{A Note on Precision} \label{S:Precision}

\noindent
Eq.\ \eqref{E:J_n-t_n} requires that we first compute the sum $s(n)$, then raise $s(n)-a_0$ to a very high power.
For example, with $a_0=7$, we have $n=168804$ and $s(168804)\approx 7+1.79178\cdot 10^{-6}$. We then compute $(s(168804)-7)^{168804}$. Many, or even {\em all}, of the significant digits of $s(n)$ will be lost if we calculate $s(n)$ to only machine precision.
Therefore, we did our calculations twice: first, we computed each $s(n)$ to 60 decimals, then used this value in Equation \eqref{E:t_n}.
Then, we repeated the calculations, this time, computing each $s(n)$ to 70 decimals.
These high-precision results agree with each other to more decimals than we show in Table \ref{Ta:SmallTable}. On the other hand, for $a_0=7$, only the first {\em three} digits of the machine precision calculation agree with these high-precision results.
Worse, when we do the calculation for $a_0=8$ in machine precision, we get, approximately,
\[
  \frac{\pi}{2} \, (1 - 1.03496 \cdot 10^{-8742942}) \, .
\]
Note that {\em all digits and the exponent} are different from the high-precision result in Table \ref{Ta:SmallTable}.

\section{Applying the Euler-Maclaurin Summation Formula} \label{S:EulerMac}

\noindent
In our special case where $a_k = 1/(2 k - 1)$ for $k \ge 1$, we can use estimates of partial sums of the harmonic series to estimate $s(n)=\sum_{k=1}^n a_k$.
Define $H_N$ to be the $N^{\text{th}}$ partial sum of the harmonic series
\[
  H_N = \sum_{k=1}^N \frac{1}{k} \, .
\]
Notice that
\begin{equation} \label{E:SumAkAsHarmonicSum}
  \sum_{k=1}^n a_k = \sum_{k=1}^n \frac{1}{2 k - 1}
   = H_{2n-1} - \frac{1}{2} H_{n-1} \, .
\end{equation}
$H_N$ has the asymptotic approximation (see \cite{WikiHarmonicNumber})
% stack exchange claims this is also proved in Chapter 9 of Concrete mathematics by Graham, Knuth, and Patashnik.
\[
  H_N \asymp \ln(N) + \gamma + \frac{1}{2N} - \sum_{k=1}^{\infty} \frac{B_{2k}}{2k N^{2k}}
  = \ln(N) + \gamma +  \frac{1}{2N} -  \frac{1}{12N^2} +  \frac{1}{120N^4} - \ldots \,
\]
which is proved (see \cite[p.\ 78]{Havil}) using the Euler-Maclaurin summation formula.
We can use this in Eq. \eqref{E:SumAkAsHarmonicSum} to get a good approximation to $s(n)=\sum_{k=1}^n a_k$.

Here, we show how to use the Euler-Maclaurin summation formula to calculate exactly the smallest $n$ for which the sum of the $a_k$ exceeds $a_0$. We also use the Euler-Maclaurin summation formula to calculate the sum of the $a_k$ for any large $n$. This method applies to general $a_k$, and gives us an estimate of the error.

One version of the Euler-Maclaurin summation formula is (see e.g. \cite[pp.\ 542-543]{Knopp})
\begin{align} \label{E:Euler-Maclaurin}
  \sum_{k=m}^n f(k)
  = {} & \int_m^n f(x)\,\dd x + \frac{f(m)+f(n)}{2}\nonumber\\
  & + \sum_{j=1}^\mu\frac{B_{2j}}{(2j)!}
  \left(f^{(2j-1)}(n)-f^{(2j-1)}(m)\right) + R_\mu(m,n)
\end{align}
with the remainder term
\begin{align} \label{E:Remainder}
  R_\mu(m,n) 
  = {} & \int_m^n\frac{B_{2\mu+1}(x-\lfloor x\rfloor)}{(2\mu+1)!}\,
  f^{(2\mu+1)}(x)\:\dd x\nonumber\\[0.1cm]
  = {} & \frac{1}{(2\mu+1)!}\,\sum_{k=m}^{n-1}\int_0^1 B_{2\mu+1}(x)\,
  f^{(2\mu+1)}(k+x)\:\dd x\,. 
\end{align}
$B_k(x)$ denotes the $k^\tn{th}$ Bernoulli polynomial, and $B_k=B_k(0)$ the $k^\tn{th}$ Bernoulli number. In our case we have
\beq
  f(x) = \frac{1}{2x-1}\qquad\mbox{and}\qquad
  a_k = f(k) = \frac{1}{2k-1}\,.  
\eeq 
Now we will derive an estimate of $R_\mu(m,n)$. For the $k^\tn{th}$ derivative of $f$, one finds
\beq
  f^{(k)}(x) = \frac{(-1)^k\,2^k\,k!}{(2x-1)^{k+1}}\,.
\eeq
Since all the functions $|f^{(k)}(x)|$, $k=0,1,2,\ldots$, are strictly decreasing, for the terms in the sum of \eqref{E:Remainder} we find
\beq
  \left| \int_0^1 B_{2\mu+1}(x) \, f^{(2\mu+1)}(k+1+x) \: \dd x \right| 
< \left| \int_0^1 B_{2\mu+1}(x )\, f^{(2\mu+1)}(k+x) \: \dd x \right| \, .
\eeq
The absolute value of each integral on the right-hand side of Equation \eqref{E:Remainder} is at most
\begin{equation*}
\left| \int_0^1 B_{2\mu+1}(x )\, f^{(2\mu+1)}(m+x) \: \dd x \right|
\end{equation*}
and there are $n-m$ of these integrals.
Therefore
\beqn \label{E:Estimation}
  \big|R_\mu(m,n)\big| < \big|\widetilde{R}_\mu(m,n)\big|
\eeqn
where
\beq
  \widetilde{R}_\mu(m,n) 
  = \frac{n-m}{(2\mu+1)!}\int_0^1 B_{2\mu+1}(x)\,f^{(2\mu+1)}(m+x)\:\dd x \, .
\eeq
Equation \eqref{E:Estimation} is the desired estimate for $R_\mu(m,n)$. Furthermore, all
\beq
  \int_0^1 B_{2\mu+1}(x)\,f^{(2\mu+1)}(k+x)\:\dd x\,,\; k=m,\ldots,n-1,
\eeq
have the same sign which is equal to the sign of $R_\mu(m,n)$, and to the sign of $\widetilde{R}_\mu(m,n)$. 

Using the integral
\beq
  \int_m^n f(x)\,\dd x
= \int_m^n\frac{\dd x}{2x-1}
= \frac{1}{2}\,[\ln(2n-1)-\ln(2m-1)]\,,
\eeq
we get the explicit summation formula
\beqn \label{E:expl_sum_form1}
 \sum_{k=m}^n\frac{1}{2k-1} = \varphi_\mu(m,n) + R_\mu(m,n)
\eeqn
with the approximation
\begin{align} \label{E:phi}
  \varphi_\mu(m,n)
= {} & \frac{1}{2}\left(\ln(2n-1)-\ln(2m-1)+\frac{1}{2m-1}+\frac{1}{2n-1}\right)\nonumber\\[0.1cm]
& - \sum_{j=1}^\mu\frac{2^{2j-1}B_{2j}}{2j}
    \left(\frac{1}{(2n-1)^{2j}}-\frac{1}{(2m-1)^{2j}}\right)    
\end{align}
and the remainder term
\beq
  R_\mu(m,n)
= -2^{2\mu+1}\,\sum_{k=m}^{n-1}
  \int_0^1\frac{B_{2\mu+1}(x)}{[2(k+x)-1]^{2\mu+2}}\,\dd x \, .
\eeq
The explicit formula for the error bound is
\beqn \label{E:error_bound}
  \widetilde{R}_\mu(m,n)
= -2^{2\mu+1}\,(n-m)
  \int_0^1\frac{B_{2\mu+1}(x)}{[2(m+x)-1]^{2\mu+2}}\,\dd x \, .
\eeqn

%\subsection{ Calculating the {\bf {\em n}} that satisfies \eqref{E:SumExceedsa0} } \label{S:CalculatingN}
\subsection{ Calculating the {\bf {\em n}} That Satisfies Inequalities (\ref{E:SumExceedsa0}) } \label{S:CalculatingN}

\noindent

Our first application of Eq.\ \eqref{E:Euler-Maclaurin} is to calculate the value of $n = n_0$ that satisfies \eqref{E:SumExceedsa0} for a fixed integer value of $a_0$.
(We can calculate the integrals for $a_0 < 10$ without too much trouble, so here, we are interested in the values $a_0 \geq 10$.)
Writing the sum $s(n)$ as
\beq
  s(n)
= \sum_{k=1}^n \frac{1}{2 k - 1}
= 1 + \frac{1}{3} + \frac{1}{5} + \cdots + \frac{1}{2 n-1}
= \frac{1}{1} + \frac{1}{1+1\cdot 2} + \frac{1}{1+2\cdot 2}
+ \cdots + \frac{1}{1+(n-1)2}\,, 
\eeq
from a theorem of Nagell \cite[pp.\ 10-14]{Nagell} (see also \cite{Belbachir}) it easily follows that $s(n)$ is never an integer except $s(1) = 1$. So we can replace \eqref{E:SumExceedsa0} by
\beqn \label{E:n-condition}
  s(n-1) < a_0 < s(n)\,.
\eeqn
Using Eq.\ \eqref{E:expl_sum_form1}, we have
\beqn \label{E:s(n)-2}
  s(n) = s(m-1)+\varphi_\mu(m,n)+R_\mu(m,n)\,.
\eeqn
Hence an approximation for $s(n)$ is
\beqn \label{E:tilde-s}
  \tilde{s}_{m,\,\mu}(n) = s(m-1)+\varphi_\mu(m,n)\,.
\eeqn
%The exact sum $s(m-1)=\sum_{k=1}^{m-1}(2k-1)^{-1}$ will be used to achieve the required precision of $\tilde{s}_{m,\,\mu}(n)$.
The error bound for the sum
\[
\sum_{k=m}^n \frac{1}{2 k - 1}
\]
is given by
\beqn \label{E:error-bound-condition}
  \big| \tilde{s}_{m, \, \mu}(n)-s(n) \big|
< \big| \widetilde{R}_\mu(m,n) \big| \, .
\eeqn
Eq.\ \eqref{E:error_bound} shows that a larger $m$ (with $n$ and $\mu$ fixed) makes the error bound smaller.
Therefore, we begin by explicitly computing the sum of the first $m-1$ terms,
\[
s(m-1) = \sum_{k=1}^{m-1} \frac{1}{2 k - 1}
\]
to high precision.
This will be used to achieve the required precision of $\tilde{s}_{m,\,\mu}(n)$. 

Replacing the integer variable $n$ in Eq.\ \eqref{E:tilde-s} by the real variable $x$, we get the equation
\beqn \label{E:tilde-s_with_x}
  \tilde{s}_{m,\,\mu}(x) = s(m-1)+\varphi_\mu(m,x)
\eeqn
with parameters $m$ and $\mu$. Solving
\beqn \label{E:Calc-n}
%  \tilde{s}_{m,\,\mu}(x) = s(m-1)+\varphi_\mu(m,x) = a_0\,,
  s(m-1) + \varphi_\mu(m,x) = a_0
\eeqn
gives a value of $x$ that approximates $n_0$.
Call this root $r$.

%\begin{figure}[h]
\begin{figure}[ht]
  \vspace{0.2cm}
  \begin{center}
    \includegraphics[scale=1]{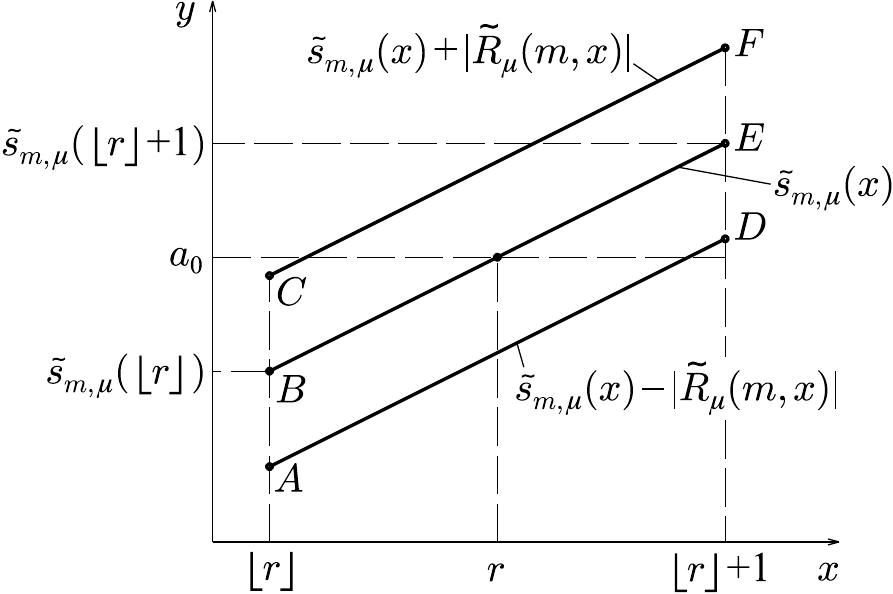}
  \end{center}
  \vspace{-0.2cm}
  \caption{\label{Fig_1} Checking if $\lfloor r \rfloor + 1 = n_0$}
\end{figure}

Now we must find a criterion that allows us to check if $\lfloor r \rfloor+1$ is the value $n_0$ that satisfies \eqref{E:n-condition}. To this purpose we consider Fig.\ \ref{Fig_1} with the graphs of the functions
\beq
  \tilde{s}_{m,\,\mu}(x)
	+ \big|\widetilde{R}_\mu(m,x)\big|\,,\qquad
  \tilde{s}_{m,\,\mu}(x)\,,\qquad
  \tilde{s}_{m,\,\mu}(x)
	- \big|\widetilde{R}_\mu(m,x)\big|\,, 
\eeq  
where (cf.\ Eq.\ \eqref{E:error_bound})
\beqn \label{E:error_bound_with_x}
  \widetilde{R}_\mu(m,x)
= -2^{2\mu+1}\,(x-m)
  \int_0^1\frac{B_{2\mu+1}(t)}{[2(m+t)-1]^{2\mu+2}}\,\dd t \, .
\eeqn

\bigskip\bigskip
Now we distinguish the following two cases:

\begin{itemize} [leftmargin=0.6cm]
 
\item[a)] $\widetilde{R}_\mu(m,x)>0$: From \eqref{E:error-bound-condition} it follows that the point $(\lfloor r \rfloor, s(\lfloor r \rfloor))$ is a point  of the line segment $\overline{BC}$, and $(\lfloor r \rfloor + 1, s(\lfloor r \rfloor + 1))$ is a point of the line segment $\overline{EF}$. If, where $\big|\overline{BC}\big|$ denotes the length of $\overline{BC}$, 
\beqn \label{E:criterion1}
  \big|\overline{BC}\big|
= \tilde{s}_{m,\,\mu}(\lfloor r \rfloor)
+ \big|\widetilde{R}_\mu(m,\lfloor r \rfloor)\big|
- \tilde{s}_{m,\,\mu}(\lfloor r \rfloor) 
< a_0-\tilde{s}_{m,\,\mu}(\lfloor r \rfloor)\,,
\eeqn
we know that $s(\lfloor r \rfloor) < a_0 < s(\lfloor r \rfloor + 1)$, hence $n_0 = \lfloor r \rfloor + 1$. Using \eqref{E:tilde-s_with_x} we write the inequality in \eqref{E:criterion1} in the final form
\beqn \label{E:criterion1a}
  a_0 - [s(m-1)+\varphi_\mu(m,\lfloor r \rfloor)]
> \big|\widetilde{R}_\mu(m,\lfloor r \rfloor)\big|\,.
\eeqn

\item[b)] $\widetilde{R}_\mu(m,x)<0$: We have $(\lfloor r \rfloor, s(\lfloor r \rfloor)) \in \overline{AB}$ and $(\lfloor r \rfloor + 1, s(\lfloor r \rfloor + 1)) \in \overline{DE}$. If 
\beqn \label{E:criterion2}
  \big|\overline{DE}\big|
= \tilde{s}_{m,\,\mu}(\lfloor r \rfloor + 1)
- \left[\tilde{s}_{m,\,\mu}(\lfloor r \rfloor + 1)
- \big|\widetilde{R}_\mu(m,\lfloor r \rfloor + 1 )\big|\right]
< \tilde{s}_{m,\,\mu}(\lfloor r \rfloor + 1) - a_0\,,
\eeqn
we know that $s(\lfloor r \rfloor) < a_0 < s(\lfloor r \rfloor + 1)$, hence $n_0 = \lfloor r \rfloor + 1$. The inequality in \eqref{E:criterion2} may be written as
\beqn \label{E:criterion2a}
  s(m-1)+\varphi_\mu(m,\lfloor r \rfloor + 1) - a_0
> \big|\widetilde{R}_\mu(m,\lfloor r \rfloor + 1)\big|\,.
\eeqn

\end{itemize}

Inequalities \eqref{E:criterion1a} and \eqref{E:criterion2a} can be combined together into 
\beqn \label{E:criterion_combined}
  \big|s(m-1)+\varphi_\mu(m,n) - a_0\big|
> \big|\widetilde{R}_\mu(m,n)\big|
\eeqn
with
\beq
  n = 
  \begin{cases}
	\lfloor r \rfloor & \text{if \, $\widetilde{R}_\mu(m,x) > 0$}\,,\\[0.1cm]
	\lfloor r \rfloor + 1& \text{if \, $\widetilde{R}_\mu(m,x) < 0$}\,.
  \end{cases}
\eeq 

As an example we will calculate $n_0$ for $a_0=10$. This allows us to check our result against Table \ref{Ta:SmallTable}.
From Table \ref{Ta:SmallTable}, it is clear that, for $a_0\ge 10$, $n_0$ is at least several million, so we choose $m=100001$ and find
\beq
 %s(m-1) = s(200000) \approx 7.0847913357763194378 \, .
  s(m-1) = s(100000) \approx 6.73821774549790928310 \, .
\eeq
%We compute this sum to 20 decimal places, but it is easy to compute more.
In this example, we compute this sum to 20 decimal places, but it is easy to compute more.
We choose $\mu = 3$ in Eq.\ \eqref{E:tilde-s}.
For larger $a_0$, or to obtain even more decimals of the sum of the $a_k$ that exceeds $a_0 = 10$, it might be necessary to use a larger $m$, to compute $s(m-1)$ to more decimals, or to use a larger value of $\mu$ (or all of the above).

We then use {\em Mathematica}'s \verb+FindRoot+ function to solve the equation
\beq
  %s(200000) + \varphi_3(200001, x) = 10
  s(100000) + \varphi_3(100001, x) = 10
\eeq
for $x$.
Expanding $\varphi_3(m, x)$, we get
\begin{align*}
  \varphi_3(m, x) 
= {} & \frac{1}{2} \left(\frac{1}{2m-1} - \ln (2m-1) + \frac{1}{2x-1} + \ln(2x-1)\right)
- \frac{1}{6} \left(\frac{1}{(2x-1)^2} - \frac{1}{(2m-1)^2}\right) \\
& - \frac{8}{63} \left(\frac{1}{(2x-1)^6} - \frac{1}{(2m-1)^6}\right) + \frac{1}{15} \left(\frac{1}{(2x-1)^4} - \frac{1}{(2m-1)^4}\right).
\end{align*}
If we substitute $m = 100001$ and combine the numeric terms together into a decimal value, we get
% -6.1030363227671701980868916
\begin{align*}
& \varphi_3(100001, x) =    \notag \\
& -6.10303632276717019809 + \frac{1}{4x-2} - \frac{1}{6(1-2x)^2} + \frac{1}{15(1-2x)^4} - \frac{8}{63(1-2x)^6} + \frac{1}{2} \ln (2x-1)
\end{align*}
So, the equation we want to solve is
\begin{align} \label{E:EquationToSolve}
         & -6.10303632276717019809 + \frac{1}{4x-2} - \frac{1}{6(1-2x)^2} + \frac{1}{15(1-2x)^4} - \frac{8}{63 (1-2x)^6} + \frac{1}{2} \ln (2x-1)    \notag \\
         & = 10 - 6.73821774549790928310 = 3.26178225450209071690 \, .
\end{align}
% end display of actual equation
We find $x = r = 68100150.0149$.
The signed error bound is found by numerical integration of Eq.\ \eqref{E:error_bound}:
\beq
  \widetilde{R}_3(100001, 68100150.0149) \approx -1.13323 \cdot10^{-39}.
\eeq
Since this error bound is less than $0$, we use Eq.\ \eqref{E:criterion_combined} with $n = \lfloor r \rfloor + 1 = 68100151$ in order to check if $n_0 = \lfloor r \rfloor + 1$. One finds that
\beq
  \big|s(100000) + \varphi_3(100001, 68100151) - 10\big|
  \approx 7.23308281312 \cdot 10^{-9} \, ,
\eeq
so the condition \eqref{E:criterion_combined} holds true, hence $n_0 = 68100151$. This confirms the $n = 68100151$ in Table~\ref{Ta:SmallTable} that was found by brute force.

%The {\em Mathematica} module \verb+getNValueAndSumForA0[ ]+ in Appendix \ref{S:MMACode-EulerNac} will compute the value of $n_0$.
%To find $n_0$ for $a_0 = 10, 11, \ldots, 25$, we use $m = 100001$.
%To compute $n_0$ for this range of $a_0$ values, it is sufficient to use $\mu = 3$.

Once we know the value of $n = n_0$ for which the sum $s(n_0)$ first exceeds $a_0 = 10$, we must compute $s(n_0)$, which is used in Equations \eqref{E:J_n-t_n} and \eqref{E:t_n} to compute the value of the $\si$ integral.
For example, with $a_0 = 10$, we find that $n_0 = 68100151$ and
\beq
  %s(n_0) = s(100000) + \varphi_3(100001, n_0) \approx 10.00000000723308281312 \, .
  s(n_0) \approx s(100000) + \varphi_3(100001, n_0) \approx 10.00000000723308281312 \, .
\eeq

Equation \eqref{E:t_n} requires that we raise the difference $s(n_0) - a_0$ to the high power $n_0$. (Note the loss of precision that occurs when we perform this subtraction). So, we may need to compute more accurate approximations $\varphi_{\mu}(m, n_0)$ using values of $\mu > 3$. Tables \ref{Ta:NValues} and \ref{Ta:Sums} below display $n_0$ and the approximate values of $s(n_0)$ for $a_0 = 10, 11, \ldots, 25$.
To obtain these values, we use $m = 100001$ and compute $s(m-1)$ to 100 decimal places, then use $\mu = 10$ to compute each $n_0$ and $\varphi_{\mu}(m, n_0)$.
The {\em Mathematica} module \verb+getNValueAndSumForA0[ ]+ in Appendix \ref{S:MMACode-EulerNac} performs these calculations.

%Note that if we compute $s(m-1)$ to 100 decimal places, then no matter how small the error estimate $\widetilde{R}_{\mu}$ is, we will never know the sum $s(n_0) = s(m-1) + \varphi_3(m, n_0)$ to more than 100 decimal places.

%Note that if we compute $s(m-1)$ to 100 decimal places, then even if the error estimate $\widetilde{R}_{\mu}$ is less than $10^{-100}$, we can \emph{not} compute the sum $s(n_0) = s(m-1) + \varphi_3(m, n_0)$ to more than 100 correct decimal places.

%Note that if we compute the initial sum $s(m-1)$ to only $D$ decimal places, then we can never compute $s(j)$ to more than $D$ correct decimal places for any $j > m-1$.
%Therefore, even if the error estimate $\widetilde{R}_{\mu}$ is less than $10^{-D}$, we can \emph{not} compute the sum $s(n_0) = s(m-1) + \varphi_3(m, n_0)$ to more than $D$ correct decimal places.

Note that if we compute the initial sum $s(m-1)$ to only $D$ decimal places, then we can never compute $s(j)$ to more than $D$ correct decimal places for any $j > m-1$, \emph{even if} the error estimate $\big| \widetilde{R}_{\mu} \big|$ is less than $10^{-D}$.

\section{Calculating the Integrals} \label{S:CalculatingIntegrals}

For a given $a_0$, we first compute the corresponding value of $n = n_0$ and the approximate value of $s(n)$, as shown in Tables \ref{Ta:NValues} and \ref{Ta:Sums}.
The next task is to compute the value of $J_n(a_0)$ using Equations \eqref{E:J_n-t_n} and \eqref{E:t_n}. The value of $(s(n)-a_0)^n$ can easily be obtained from the approximate value of $s(n)$ in Table \ref{Ta:Sums}, although for large $n$, we must use logarithms to prevent underflow.

\begin{table}[ht]
\begin{center}
\begin{tabular}{|r|r|} \hline
\rule{0pt}{12pt}
$a_0$ & \multicolumn{1}{c|}{$n$} \\[2pt] \hline\hline
10 & 68100151 \\
11 & 503195829 \\
12 & 3718142208 \\
13 & 27473561358 \\
14 & 203003686106 \\
15 & 1500005624924 \\
16 & 11083625711271 \\
17 & 81897532160125 \\
18 & 605145459495141 \\
19 & 4471453748222757 \\
20 & 33039822589391676 \\
21 & 244133102611731231 \\
22 & 1803913190804074904 \\
23 & 13329215764452299411 \\
24 & 98490323038288832267 \\
25 & 727750522131718025058 \\ \hline
\end{tabular}
\vspace{0.2cm}
\renewcommand{\arraystretch}{1}
\caption{Values of $n = n_0$, for each $a_0$}
\label{Ta:NValues}
\end{center}
\end{table}

\begin{center}
{\Small
\begin{table}[ht]
\renewcommand{\arraystretch}{1.2}
\begin{tabular}{|c|l|l|}\hline
\rule{0pt}{16pt}
\normalsize{$a_0$} & \multicolumn{1}{c|}{\normalsize{$\tilde{s}_{100001,\,10}(n)$}} & \multicolumn{1}{c|}{\normalsize{$\widetilde{R}_{10}(100001,n)$}}\\[3pt] \hline\hline
10               &
$10 + 7.23308281311740815495440938881892875629793229610802275303838659 \cdot 10^{-9}$ & \quad $2.10 \cdot 10^{-104}$\\
11              & 
$11 + 1.93429694721571938243592220609208607459386666993511996115170447 \cdot 10^{-10}$ & \quad $1.56 \cdot 10^{-103}$\\
12             & 
$12 + 2.81704757017003986061562163359221047582420212335754506062428273 \cdot 10^{-11}$ & \quad $1.15 \cdot 10^{-102}$\\
13            & 
$13 + 1.51784528343340657974855459172890208869659078869207257172632024 \cdot 10^{-11}$ & \quad $8.50 \cdot 10^{-102}$\\
14           & 
$14 + 1.20004359101609629122080445652372955218041261117217502496026183 \cdot 10^{-12}$ & \quad $6.28 \cdot 10^{-101}$\\ 
15          &
$15 + 2.19180272149887470909606761989402891098093034885435056479463010 \cdot 10^{-13}$ & \quad $4.64 \cdot 10^{-100}$\\
16         & 
$16 + 4.03776701092542650935062088404145888641302878626731173386452500 \cdot 10^{-14}$ & \quad $3.43 \cdot 10^{-99}$\\
17         & 
$17 + 3.96811610610919880012621568968292119992007054389850812439945895 \cdot 10^{-16}$ & \quad $2.54 \cdot 10^{-98}$\\
18        & 
$18 + 6.44184629552359167120616511513071089671769035954843683552410240 \cdot 10^{-16}$ & \quad $1.87 \cdot 10^{-97}$\\
19       & 
$19 + 4.40658835470283585071853820223887285629968834008507731684755191 \cdot 10^{-17}$ & \quad $1.38 \cdot 10^{-96}$\\
20      &
$20 + 2.90258683104894913499203070153167600669577064916856020926287204 \cdot 10^{-18}$ & \quad $1.02 \cdot 10^{-95}$\\
21     & 
$21 + 7.34280669057054306832818424563959102068261955548079016234613646 \cdot 10^{-19}$ & \quad $7.56 \cdot 10^{-95}$\\
22    & 
$22 + 1.30683560567708459204537388912129731492458474471888171662329379 \cdot 10^{-19}$ & \quad $5.58 \cdot 10^{-94}$\\
23   &
$23 + 3.40633844408955109014083203224199839911999656758748815953125439 \cdot 10^{-20}$ & \quad $4.13 \cdot 10^{-93}$\\
24   & 
$24 + 3.79499658486046318316555581259771062170888781423675014763472623 \cdot 10^{-21}$ & \quad $3.05 \cdot 10^{-92}$\\
25  & 
$25 + 1.14325480646582051223669818246654129326458197624224895807362049 \cdot 10^{-22}$ & \quad $2.25 \cdot 10^{-91}$\\
\hline
\end{tabular}
\vspace{0.2cm}
\renewcommand{\arraystretch}{1}
\caption{$a_1 + a_2 + \ldots + a_n$ ($n$ from Table \ref{Ta:NValues}), and error bound from \eqref{E:error_bound}, for each $a_0$}
\label{Ta:Sums}
\end{table}
} % end small
\end{center}

Equation \eqref{E:t_n} requires that we compute the product of $a_k$ for $k = 1, \ldots, n$.
In our examples, we have $a_k = 1/(2k - 1)$ for $k = 1, \ldots, n$, so we can rewrite this product as follows:
\begin{equation*}
  \prod_{k=1}^n a_k
= \prod_{k=1}^n \frac{1}{2k-1}
= \frac{1}{2n-1} \cdot \frac{1}{2n-3} \cdots \frac{1}{3} \cdot \frac{1}{1}\,.
\end{equation*}
We can express this in terms of factorials:
\begin{align*}
  & \frac{1}{2n-1} \cdot \frac{1}{2n-3} \cdots \frac{1}{3} \cdot \frac{1}{1}
= \frac{ 2n \cdot (2n-2) \cdots 4 \cdot 2 }{ 2n \cdot (2n-1) \cdot 2 \cdot 1 }
= \frac{ 2^n \cdot n \cdot (n-1) \cdots \cdot 2 \cdot 1 }{ 2n \cdot (2n-1) \cdot 2 \cdot 1 }
= \frac{ 2^n n! }{ (2n)! } \, .
\end{align*}

%\begin{align*}
%  (2n-1)!!
%= {} & (2n-1)(2n-3)\cdots 3\cdot 1
%= \frac{2n\,(2n-1)\cdots 2\cdot 1}{2n\,(2n-2)\cdots 4\cdot 2}\\
%= {} & \frac{2n\,(2n-1)\cdots 2\cdot 1}{2^n\,n\,(n-1)\cdots 2\cdot 1}
%= \frac{(2n)!}{2^n\,n!}\,.  
%\end{align*}
%Therefore, Eq.\ \eqref{E:J_n-t_n-2} can be written as
Therefore, we can rewrite Equations \eqref{E:J_n-t_n} and \eqref{E:t_n} as
\beqn \label{E:J_n-3}
  J_n(a_0) = \frac{\pi}{2}\,\{1-t_n(a_0)\}\qquad\mbox{where}\qquad 
%  t_n(a_0) = \frac{(s(n)-a_0)^n \, (2n)!}{2^{2n-1} \, n!^2}
  t_n(a_0) = \frac{ (s(n)-a_0)^n }{ 2^{2n-1} } \cdot \frac{ (2n)! }{ n!^2 }
\eeqn

When $n$ is as large as some of the values in Table \ref{Ta:NValues}, we will need to use logarithms. Taking natural logarithms, we have
\beqn \label{E:lnTiny}
  \ln t_n(a_0) = n\ln (s(n)-a_0)-(2n-1)\ln 2+\ln((2n)!)-2\ln(n!)\,.
\eeqn
We denote by lg the the log base 10 and have
\beqn \label{E:log10Tiny}
  %\lg t_n(a_0) = \lg\ee\cdot\ln t_n(a_0)
  \lg t_n(a_0) = \frac{\ln t_n(a_0)}{\ln(10)}
\eeqn
with $t_n(a_0)$ from Eq.\ \eqref{E:lnTiny}.
% note: later, we show how to extract the mantissa and exponent.
% see equations {E:GetExponent} and {E:GetMantissa}
%so, this is removed from here.
%After calculating $\lg t_n(a_0)=:\ell$, we can extract the mantissa $m$ and the exponent $p$ of $t_n(a_0)$, and then display $t_n(a_0)$ in scientific notation using
%\beq
%  m = 10^{\ell - \lfloor\ell\rfloor}\,,\qquad p = \lfloor\ell\rfloor\,.
%\eeq
%Now, it remains to estimate the last two terms of Eq.\ \eqref{E:lnTiny}.
Now, it remains to estimate the logarithms of the factorials in Eq.\ \eqref{E:lnTiny}.
We will discuss two techniques: approximations based on Stirling's formula, and another application of the Euler-Maclaurin summation formula.

\subsection{Estimating {\bf {\em n}}! With Stirling's Formula}

\noindent

For large $n$, we can use Stirling's formula to get reasonable approximations to the above combination of factorials.
These approximations mainly require a few exponentiations, which is much faster than doing $O(n)$ multiplications, especially when $n > 10^6$.

The simple version of Stirling's formula is
\[
  n! \approx \left( \frac{n}{\ee} \right)^n \cdot \sqrt{2 \pi n} \, .
\]
Strictly speaking, this is not an approximation to $n!$, but is an asymptotic relation. This means that, as $n$ approaches $\infty$, the \emph{ratio} of the right-hand side to $n!$ approaches 1. However, for the large $n$ in Table \ref{Ta:Sums}, the ratio is 1 to several decimal places, so we can obtain several decimal places of the right side in $\pi[1-t_n(a_0)]/2$.

A more rigorous approach is to compute both lower and upper bounds for $n!$, based on refined versions of Stirling's formula.
These bounds were proved in \cite{Robbins}; see also \cite{MathWorldStirling}:
\begin{equation*}
  \left( \frac{n}{\ee} \right)^n \cdot \sqrt{2 \pi n} \cdot \ee^{\frac{1}{12n + 1}}
  < n! <
  \left( \frac{n}{\ee} \right)^n \cdot \sqrt{2 \pi n} \cdot \ee^{\frac{1}{12n}} \, .
\end{equation*}

The following improved lower bound was proved in \cite{Maria}, so we will use
\begin{equation} \label{E:FactorialBounds}
  \left( \frac{n}{\ee} \right)^n \cdot \sqrt{2 \pi n} \cdot \ee^{\frac{1}{12n + \frac{3}{4n + 2}}}
  < n! <
  \left( \frac{n}{\ee} \right)^n \cdot \sqrt{2 \pi n} \cdot \ee^{\frac{1}{12n}} \, .
\end{equation}

How close are the lower and upper bounds in \eqref{E:FactorialBounds}? The ratio of the upper bound to the lower bound is
\begin{equation*}
  r_n
  = \exp \left( \frac{1}{12n} - \frac{1}{12n + \frac{3}{4n + 2}} \right)
  = 1 + \frac{1}{192 n^3} - \frac{1}{384 n^4} + O\left(\frac{1}{n^5}\right) \, .
\end{equation*}
For large $n$, this ratio is quite small.

For $n = 10^6$, the lower and upper bounds in \eqref{E:FactorialBounds} are
\begin{align*}
  & 8.2639316883312400623566 \cdot 10^{5565708} \text{ and} \\
  & 8.2639316883312400623996 \cdot 10^{5565708}
\end{align*}
Notice that the first 20 digits of the lower and upper bounds are the same.
Therefore, we know that $n!$ begins these 20 digits.
In fact, \textit{Mathematica} can calculate $n!$, and the value is about $8.2639316883312400623766 \ldots \cdot 10^{5565708}$.

(But note that if the lower and upper bounds of some $x$ were 1.9 and 2.1, then \emph{no} digits would agree, but we would still know that $x = 2 \pm 0.1$.)

We conclude that, using the bounds in \eqref{E:FactorialBounds}, we can compute at least the first 10 significant digits of factorials of very large numbers.

%When $n$ is as large as some of the values in Table \ref{Ta:Sums}, we will need to use logarithms. %already said this above!
Although the natural logarithms of the bounds above would be a little simpler, logs base 10 will most easily produce displayable values.
Taking logs in \eqref{E:FactorialBounds}, we have
\begin{equation*}
    n \lg \left( \frac{n}{\ee} \right) + \frac{ \lg( 2 \pi n ) }{2} + \frac{\lg(\ee)}{12n + \frac{3}{4n + 2}}
  < \lg( n! )
  < n \lg \left( \frac{n}{\ee} \right) + \frac{ \lg( 2 \pi n ) }{2} + \frac{\lg(\ee)}{12n} \, .
\end{equation*}
So, for a given $n$, define $b_1(n)$ to be the lower bound of $\lg(n!)$:
\begin{equation*}
  b_1(n) = n \lg \left( \frac{n}{\ee} \right) + \frac{ \lg( 2 \pi n ) }{2} + \frac{\lg(\ee)}{12n + \frac{3}{4n + 2}} \, ,
\end{equation*}
and define $b_2(n)$, as the upper bound of $\lg(n!)$:
\begin{equation*}
  b_2(n) = n \lg \left( \frac{n}{\ee} \right) + \frac{ \lg( 2 \pi n ) }{2} + \frac{\lg(\ee)}{12n} \, .
\end{equation*}

Our goal is to compute $\lg((n!)^2/(2n)!)$.
To get a lower bound, we use the lower bound for $\lg(n!)$ and the upper bound for $\lg((2n)!)$.
So, the lower bound for $\lg((n!)^2/(2n)!)$ is
\begin{equation} \label{E:LogLowerBoundRatio}
  2 b_1(n) - b_2(2n) \, ,
\end{equation}
and the respective upper bound is
\begin{equation} \label{E:LogUpperBoundRatio}
  2 b_2(n) - b_1(2n) \, .
\end{equation}

Now suppose we have computed $z \approx \lg(N)$ for some large $N$ (for example, $N = n!$).
To display the approximate value of $N = 10^z$, we need not compute $10^z$, which would overflow if $N$ is large enough.
Instead, we can extract the mantissa and the exponent of $z$, and then display $N$ in scientific notation.
We have
\begin{align}
  p & = \lfloor z \rfloor  \label{E:GetExponent} \\
  m & = 10^{ z - p }  \label{E:GetMantissa} \, .
\end{align}
We can now display $N$, the antilog of $z$, as $N = m \cdot 10^p$, where $1 \leq m < 10$.
However, note that if the integer exponent $p$ has $d$ significant digits, then the mantissa $m$ will have about $d$ fewer significant digits than $z$ has.
This happens because the subtraction in \eqref{E:GetMantissa} causes a loss of precision.
For example, if $z = 100000.12$, then $10^z = 10^{0.12} \cdot 10^{100000} \approx 1.3 \cdot 10^{100000}$.
Although $10^{0.12} \approx 1.3182567 \dots$, only the first two digits in the mantissa of $10^z$ are meaningful.
This is because all we know about this ``0.12'' is that it is a number between .115 and .125, and $10^{.115} \approx 1.303$ and $10^{.125} \approx 1.334$.

If $n = 10^6$, then we can estimate $\lg((n!)^2/(2n)!)$ as follows.
The lower and upper bounds for $\lg$ of $(n!)^2/(2n)!$ in \eqref{E:LogLowerBoundRatio} and \eqref{E:LogUpperBoundRatio} are
\begin{align*}
&-602056.74275297175655031271186 \quad \text{ and } \\
&-602056.74275297175655031270705 \, .
\end{align*}
These agree to 25 digits.
If we extract the mantissa and the exponent for these two logarithms, we obtain, respectively,
\begin{align*}
& 1.8082023454706427717249 \cdot 10^{-602057} \quad \text{ and } \\
& 1.8082023454706427717449 \cdot 10^{-602057} \, .
\end{align*}
These agree with each other to 20 digits.
This implies that, to 20 significant digits, the value of of $(n!)^2/(2n)!$ is $1.8082023454706427717 \cdot 10^{-602057}$.
Notice that we lost several significant digits because the exponent has 6 digits.

$n = 10^6$ is small enough that \textit{Mathematica} can compute $(n!)^2/(2n)!$ directly.
The value is about $1.8082023454706427717343 \cdot 10^{-602057}$.

Here is a straightforward implementation of some of the above equations as \textit{Mathematica} code:
{
\SMALL
\begin{verbatim}
lowerFactorial[n_] := (n/E)^n * Sqrt[2 Pi n] * Exp[1/(12 n + 3/(4 n + 2))]
upperFactorial[n_] := (n/E)^n * Sqrt[2 Pi n] * Exp[1/(12 n)]
log10LowerFactorial[n_] := n * Log[10, n/E] + Log[10, 2 Pi n]/2 + Log[10, E]/(12 n + 3/(4 n + 2))
log10UpperFactorial[n_] := n * Log[10, n/E] + Log[10, 2 Pi n]/2 + Log[10, E]/(12 n)
(* here are the lower and upper bounds of log10[ (n!)^2/(2n)! ] *)
log10LowerRatio[n_] := 2*log10LowerFactorial[n] - log10UpperFactorial[2 n]
log10UpperRatio[n_] := 2*log10UpperFactorial[n] - log10LowerFactorial[2 n]
(* log of lower and upper bounds of tiny *)
log10TinyLower[a0_, n_, s_] := n * Log[10, s - a0] - (2 n - 1) * Log[10, 2] - log10UpperRatio[n]
log10TinyUpper[a0_, n_, s_] := n * Log[10, s - a0] - (2 n - 1) * Log[10, 2] - log10LowerRatio[n]
getME1[x_] := { 10^(x - Floor[x]) , Floor[x] } (* get matissa and exponent of antilog base 10 *)
\end{verbatim}
}  % end small

Let's use this code to calculate $J_{68100151}(10)$. First, obtain $n = 68100151$ and the sum
\[
%  s = 10 + 7.2330828131174081549544093888189287562 \cdot 10^{-9}
%  s = 10 + 7.2330828131174081549544093888189287562979322961080227530383865861205404748019681208202156369 \cdot 10^{-9}
   s = 10 + 7.233082813117408154954409388818928756297 \cdot 10^{-9}
\]
truncated from Table \ref{Ta:Sums}.
(Or, one may compute $s$ by adding 68100151 terms directly).
Then, run the following \textit{Mathematica} code:
\begin{verbatim}
a0 = 10
n = 68100151
s = 10 + 7.233082813117408154954409388818928756297 * 10^(-9)
logt1 = log10TinyLower[a0, n, s]
logt2 = log10TinyUpper[a0, n, s]
{m1, e1} = getME1[logt1]
{m2, e2} = getME1[logt2]
\end{verbatim}
The results are
\begin{verbatim}
  {m1, e1} = {9.6492736004286844634795529419197, -554381308}
  {m2, e2} = {9.6492736004286844634795532800687, -554381308}
\end{verbatim}
Of the 32 significant digits in the mantissas, 24 of them agree.
Therefore, we know that
\[
t \approx 9.64927360042868446347955 \cdot 10^{-554381308} \, ,
\]
where all 24 displayed digits in the mantissa are correct.
This is consistent with the result
\[
t \approx 9.649273600428684463479553 \ldots \cdot 10^{-554381308} \, ,
\]
given in Table \ref{Ta:SmallTable}.

Finally, the following \textit{Mathematica} code will calculate $t$ for $a_0 = 25$:
% we no longer have this many digits in the sum in the wide table; use the less-accurate value in the verbatim below
%s = 25 + 1.1432548064658205122366981824665412932645819762422%48958073620493072803557442617 * 10^(-22)
%s = 25 + 1.1432548064658205122366981824665412932645819762422489580736204930728 * 10^(-22)
\begin{verbatim}
a0 = 25
n = 727750522131718025058
s = 25 + 1.1432548064658205122366981824665412932645\
8197624224895807362049 * 10^(-22)
logt1 = log10TinyLower[a0, n, s]
logt2 = log10TinyUpper[a0, n, s]
{m1, e1} = getME1[logt1]
{m2, e2} = getME1[logt2]
\end{verbatim}

Both pairs \verb+{m1, e1}+ and \verb+{m2, e2}+ are
%\begin{verbatim}
%  {2.72384864752823356168556312789934977164693139593229092635,
%   -15968197862152240928105}
%\end{verbatim}
\begin{verbatim}
  {2.7238486475282335616855631278993497716469,
   -15968197862152240928105}
\end{verbatim}
which means that
\[
  %t \approx 2.72384864752823356168556312789934977164693139593229092635 \cdot 10^{-15968197862152240928105}
   t \approx 2.7238486475282335616855631278993497716469 \cdot 10^{-15968197862152240928105}
\]
where all 41 digits in the mantissa are correct.
% we can confirm that all of these digits are correct. here are 64 digits:
% for a0 = 25, we have t = 2.723848647528233561685563127899349771646931395932290926352454516 * 10^-15968197862152240928105,
% which can be calculated using
% m1 = 300001, mu1 = 14, nDecimals1 = 110, accGoal1 = 25, workPrec1 = 50
% m2 = 200001, mu2 = 7, nDecimals2 = 110, workPrec2 = 45

The number of correct digits here is determined by how closely the lower and upper bounds of $n!$ agree.
%As $n$ increases, so does the number of digits to which the lower and upper bounds agree.
%The limitation here is that we cannot arbitrarily increase the accuracy of these approximations.
The limitation here is that we are stuck with whatever accuracy these approximations provide.
The Euler-Maclaurin summation formula, described next, allows us to get as much accuracy as we want.

\subsection{Estimating {\bf {\em n}}! With the Euler-Maclaurin Summation Formula}

\noindent

We can write
\beq
  \sigma(n) := \ln((2n)!)-2\ln(n!) = \sum_{k=1}^{2n}\ln k-2\,\sum_{k=1}^n\ln k\,.
\eeq
To get a good estimate for $\sigma(n)$, we will use the exact sum of $m-1$ initial terms. Therefore, we split the sums: 
\beqn \label{E:sigma}
  \sigma(n)
= \left[\sum_{k=1}^{m-1}+\sum_{k=m}^{2n}-\:2\left(\sum_{k=1}^{m-1}+\sum_{k=m}^n\right)\right]\ln k
= -\left(\sum_{k=1}^{m-1}+\sum_{k=m}^n-\sum_{k=n+1}^{2n}\right)\ln k\,.
\eeqn
It remains to estimate $\sum_{k=m}^n\ln k$ and $\sum_{k=n+1}^{2n}\ln k$. Therefore, we apply the Euler-Maclaurin summation formula \eqref{E:Euler-Maclaurin} with
\beq
  f(x) = \ln x\,.
\eeq
The derivatives are given by
\beqn \label{E:Derivatives_of_ln}
  f^{(k)}(x) = \frac{(-1)^{k-1}\,(k-1)!}{x^k}\,,\quad 
  k = 1,2,\ldots\,.
\eeqn
Furthermore, we have
\beq
  \int_m^n f(x)\,\dd x
= \int_m^n\ln x\,\dd x = x(\ln x-1)\,\Big|_m^n
= n(\ln n-1)-m(\ln m-1)\,.
\eeq
This yields
\beq
 \sum_{k=m}^n\ln k = \psi_\mu(m,n) + R_\mu^*(m,n)
\eeq
with the approximation
%\beq
%  \psi_\mu(m,n) = n(\ln n-1) - m(\ln m-1) + \frac{\ln m+\ln n}{2}\\
%  + \sum_{j=1}^\mu\frac{B_{2j}}{2j(2j-1)}
%   \left(\frac{1}{n^{2j-1}}-\frac{1}{m^{2j-1}}\right)
%\eeq
\begin{align} \label{E:factorialPsi}
  \psi_\mu(m,n) = & n(\ln n-1) - m(\ln m-1) + \frac{\ln m+\ln n}{2}  \notag \\
  & + \sum_{j=1}^\mu\frac{B_{2j}}{2j(2j-1)}
   \left(\frac{1}{n^{2j-1}}-\frac{1}{m^{2j-1}}\right)
\end{align}
and the remainder term
\begin{equation*}
  R_\mu^*(m,n)
= \frac{1}{2\mu+1}\,\sum_{k=m}^{n-1}\int_0^1\frac{B_{2\mu+1}(x)}{(k+x)^{2\mu+1}}\,\dd x\,.
\end{equation*}
Since the absolute values of all derivatives in Equation \eqref{E:Derivatives_of_ln} are strictly decreasing, we find the error estimate 
\begin{equation*} %\label{E:Estimation2}
  \big|R_\mu^*(m,n)\big| < \big|\widetilde{R}_\mu^*(m,n)\big|
\end{equation*}
with
\begin{equation} \label{E:factorialErrorR}
  \widetilde{R}_\mu^*(m,n)
= \frac{n-m}{2\mu+1}\,\int_0^1\frac{B_{2\mu+1}(x)}{(m+x)^{2\mu+1}}\,\dd x\,.
\end{equation}
It follows that the approximation for \eqref{E:sigma} is given by
\begin{equation} \label{E:factorialRatio}
  \tilde{\sigma}_{m,\,\mu}(n) = -\sum_{k=1}^{m-1}\ln k - \psi_\mu(m,n) + \psi_\mu(n+1, 2n) 
\end{equation}
with the error bound
\begin{equation} \label{E:factorialRatioErrorBound}
 \big|\tilde{\sigma}_{m,\,\mu}(n)-\sigma(n)\big|
< \big|\widetilde{R}_\mu^*(m,n)\big| + \big|\widetilde{R}_\mu^*(n+1,2n)\big|\,.
\end{equation}

%\subsection{Computing $t$}
\subsection{Computing {\bf {\em t}}}

\noindent

We can now put all of this together to compute the value of $t$.

The \textit{Mathematica} code to perform these calculations is in Appendix \ref{S:MMACode-EulerNac}.

There are three main parts of the code.

The first part, in module \verb+getNValueAndSumForA0[ ]+, computes the smallest $n = n_0$ for which the sum of $a_1 + ... + a_n$ exceeds $a_0$.
The value of the sum is also computed and saved.

The second part, in module \verb+lnFactRatio[ ]+, computes the natural logarithm of $(2n)!/(n!)^2$ for the $n = n_0$ just we found.

The third part runs a loop over some range of values of $a_0$.
This loop calls the above modules and saves and prints the results.

Appendix \ref{S:MMACode-Results} displays the results of running the \textit{Mathematica} code.

\vspace{0.25cm}

\noindent
{\bf Uwe B\"asel}, HTWK Leipzig, Faculty of Mechanical and Energy Engineering, Leipzig, Germany, uwe.baesel@htwk-leipzig.de\\[0.3cm]
{\bf Robert Baillie}, State College, Pennsylvania, USA, rjbaillie@frii.com 

\vspace{0.25cm}

\pagestyle{empty} % stop page numbering here to keep page numbers out of the mathematica code
%\newpage
\appendix
\section{\emph{Mathematica} Code} \label{S:MMACode}

%This code works \emph{at least} in \textit{Mathematica} versions 7, 8, and 9.
This code has been tested, and works, in \textit{Mathematica} versions 7, 8, and 9.

%\subsection{The Euler-Maclaurin Summation Formula} \label{S:MMACode-EulerNac}
\subsection{Code} \label{S:MMACode-EulerNac}

%\noindent

Here is the \textit{Mathematica} code to compute the integrals for $a_0 = 10$ through $a_0 = 25$.

This code uses the Euler-Maclaurin summation formula to compute the corresponding $n$ for which
\[
  \sum_{k=1}^n a_k = \sum_{k=1}^n \frac{1}{2 k - 1}
\]
first exceeds $a_0$.
The value of this sum is also computed.

For each $a_0$ and the corresponding (large) $n$, this code then uses the Euler-Maclaurin summation formula to compute logarithm of the ratio of factorials
\[
\frac{ (2 n)! }{ n!^2 } \, ,
\]
which occurs in Equation \eqref{E:J_n-3}.

This code is presented in a format (for example, without page numbers) that the user can copy and paste directly into \textit{Mathematica}.

\SMALL
\begin{lstlisting}

(*
for a given a0, the integral for a0 is
  (Pi/2) * (1 - t).
the code below computes the value of t, which will be very tiny.

the number of significant digits in t depends on the values of
m1, mu1, nDecimals1, m2, mu2, nDecimals2, and accGoal and workPrec.
for convenience, we set all of those values here.
*)


(* these are used to compute sums of ak, and the n value, given a0 *)
mu1 = 10;  (* number of derivative terms to find n and the sum of ak *)
m1 = 100001;  (* 1 + number of initial terms in sum of ak *)
nDecimals1 = 100;  (* get the sum if ak to this many digits after the decimal point *)

(* accGoal and workPrec help get more accurate roots *)
accGoal1 = 20;
workPrec1 = 2*accGoal1;


(* these are used to compute factorials of large numbers *)
mu2 = 5;   (* number of derivative terms to compute factorials *)
m2 = 100001;  (* 1 + number of initial terms in sum of logs *)
nDecimals2 = 100;  (* want the sum to this many digits after the decimal point *)

(* workPrec2 helps get a more accurate value for the integral in R2 *)
workPrec2 = 40;


(* the number of accurate digits in the result depends on these initial values *)
Print["m1 = ", m1, ", mu1 = ", mu1, ", nDecimals1 = ", nDecimals1,
  ", accGoal1 = ", accGoal1, ", workPrec1 = ", workPrec1];
Print["m2 = ", m2, ", mu2 = ", mu2, ", nDecimals2 = ", nDecimals2, ", workPrec2 = ", workPrec2];




(* define two utility functions, getME[ ] and removeQuestionableDigits[ ] *)
getME[c_] :=
Module[
  (* get matissa and exponent of c, an antilog base 10. *)
  { expo, diff, mant },
  expo = Floor[c];
  diff = c - expo;
  If[Accuracy[diff] < 2,
    mant = 1,  (* not enough significant digits remain after subtraction *)
    mant = 10^diff
  ];
  Return[ { mant , expo } ]
]  (* end of Module *)


removeQuestionableDigits[c_, errorEst_] :=
Module[
(* c may have many digits that are significant as far as Mathematica is concerned,
   but the estimated error for c, based on an integral, might make some of those 
   digits be meaningless.
   example: if
     c = 1.2345678901234567890123456789012345 * 10^7 (35 digits),
   and the error estimate for c is 6.513 * 10^-20, then return
     c = 1.2345678901234567890123457*10^7 (26 digits).
   to round c to have n digits after the decimal point, call
     removeQuestionableDigits[ c , 10^-(n+1) ] .
   note that this may introduce an error of up to 0.5 * 10^-(n+1) .
*)
  { digitsRightOfDP, log10c, mant, expo, numDigitsCorrect, power10ToRound },
  digitsRightOfDP = Floor[ -Log[10, Abs[errorEst]] ];
  If[digitsRightOfDP >= Floor[Accuracy[c]],
    Return[c]  (* nothing to do; example: c = 1.2345, errorEst = 10^-20 *)
  ];
  log10c = Log[10, Abs[c]];
  {mant, expo} = getME[log10c];  (* c = mant * 10^expo *)
  numDigitsCorrect = expo + digitsRightOfDP;
  numDigitsCorrect = Max[numDigitsCorrect, 1];
  power10ToRound = N[10^-numDigitsCorrect, numDigitsCorrect];
  Return[ 10^expo * Round[mant, power10ToRound ] ]
]  (* end of Module *)




(* here is the code to compute n for a given value of a0 *)

(* compute s(m1 - 1); see Equation (4) *)
initialAkSum = N[Sum[1/(2 k - 1), {k, 1, m1 - 1}], nDecimals1 + 10];
(*
we now have the sum to at least nDecimals1 digits, with essentially no roundoff error.
next, round it to have (nDecimals1) digits after the decimal point.
note that this sum might now have an error of 0.5*10^-(nDecimals1 + 1).
*)
initialAkSum = removeQuestionableDigits[initialAkSum, 10^-(nDecimals1 + 1)];
Print["sum of the first ", m1 - 1, " ak values = ", initialAkSum];


(* curlyPhi = Euler-Maclaurin sum without the error term; see Equation (23) *)
curlyPhi[mu_, m_, n_] := 
  1/2 (Log[2 n - 1] - Log[2 m - 1] + 1/(2 m - 1) + 1/(2 n - 1)) - 
   Sum[(2^(2 j - 1) BernoulliB[2 j])/(2 j) * (1/(2 n - 1)^(2 j) - 1/(2 m - 1)^(2 j)),
       {j, 1, mu}];


(* R from Equation (24) *)
R1[mu_, m_, n_] := -2^(2 mu + 1) (n - m) *
 ( NIntegrate[BernoulliB[2 mu + 1, x]/(2 (m + x) - 1)^(2 mu + 2), {x, 0, 1/2},
     WorkingPrecision -> workPrec1] +
   NIntegrate[BernoulliB[2 mu + 1, x]/(2 (m + x) - 1)^(2 mu + 2), {x, 1/2, 1},
     WorkingPrecision -> workPrec1] );


getNValueAndSumForA0[a0_, m_, mu_, initialAkSum_, accGoal_, workPrec_, iPrint_] :=
Module[
  (* this module computes and returns three values:
     (1) the smallest n value for which the sum of 1/ak (k >= 1) first exceeds a0;
     (2) the corresponding sum; this is slightly larger than a0;
     (3) the absolute value of error estimate of the sum.
     the error bound is based on equation (24). based on this error bound,
     we discard any digits in the sum that might not be correct.
     if iPrint == 1, this prints out various internal values for debugging.
  *)
  { rootGuess = 10^6, root, r, nValue, errorEst, sum1Est, sum1Shortened,
    sum2Est, sum2Shortened },

  sum1Est = sum1Shortened = sum2Est = sum2Shortened = 0;

  root = n /. FindRoot[initialAkSum + curlyPhi[mu, m, n] == a0, {n, rootGuess},
    AccuracyGoal -> accGoal, WorkingPrecision -> workPrec];

  (* check if the difference in r is greater than the error bound *)
  If[R1[mu, m, root] > 0,
    r = Floor[root],
    r = Floor[root] + 1
  ];
  errorEst = Abs[R1[mu, m, r]];
  If[Abs[a0 - (initialAkSum + curlyPhi[mu, m, r])] > errorEst,
    nValue = Floor[root] + 1,
    nValue = 0  (* failed to find a valid root *)
  ];

  If[iPrint == 1,
    (* R1[mu, m, r] and R1[mu, m, root] usually agree to about 9 significant digits *)
    Print["a0 = ", a0, ", root = ", root ", nValue = ", nValue,
          ", R1 = ", N[R1[mu, m, root], 6],
          ", diff1 = ", N[a0 - (initialAkSum + curlyPhi[mu, m, r]), 10],
          ", diff2 = ", N[initialAkSum + curlyPhi[mu, m, r + 1] - a0, 10]
    ];
  ];

  If[nValue > 0,  (* we do not use sum1Est or sum1Shortened *)
    (* sum1Est = initialAkSum + curlyPhi[mu, m, nValue - 1]; *)  (* a1 + a2 + ... + a[nValue-1] *)
    (* note: we know initialAkSum to only (nDecimals1) decimal places.
       therefore, we cannot know sum2Est to more decimal places than that.
       when Mathematica adds curlyPhi to initialAkSum, it will not retain
       more than (nDecimals1) decimal places in the result.
    *)
    sum2Est = initialAkSum + curlyPhi[mu, m, nValue];  (* a1 + a2 + ... + a[nValue] > a0 *)
    (* keep only those digits that are justified, given the error estimate from R1 *)
    (* sum1Shortened = removeQuestionableDigits[sum1Est, errorEst]; *)
    sum2Shortened = removeQuestionableDigits[sum2Est, errorEst];
    If[iPrint == 1,
      Print["a0 = ", a0, ", errorEst = ", N[errorEst, 6],
            ", curlyPhi = ", curlyPhi[mu, m, nValue],
       ", sum2Est = ", sum2Est, ", shortened sum2Est = ", sum2Shortened];
    ]
  ];

 Return[ { nValue , sum2Shortened , N[errorEst, 6] } ]

] (* end Module *)




(*
here is the code to use the Euler-Maclaurin summation formula to compute factorials.
*)


(* this sum will have about 7 digits before the decimal point, so add 15 here *)
initialLogSum = N[Sum[Log[k], {k, 1, m2 - 1}], nDecimals2 + 15];
(*
we now have the sum to at least nDecimals2 digits, with essentially no roundoff error.
next, round it to have (nDecimals2) digits after the decimal point.
note that this sum might now have an error of 0.5*10^-(nDecimals2 + 1).
*)
initialLogSum = removeQuestionableDigits[initialLogSum, 10^-(nDecimals2 + 1)];
Print["sum of the first ", m2 - 1," logs = ", initialLogSum];


(* curlyPsi = Euler-Maclaurin sum without the error term, see Equation (48) *)
curlyPsi[mu_, m_, n_] := 
  n*(Log[n] - 1) - m*(Log[m] - 1) + (Log[m] + Log[n])/2 +
   Sum[BernoulliB[2 j]/(2 j (2 j - 1)) * (1/n^(2 j - 1) - 1/m^(2 j - 1)),
       {j, 1, mu}];


(* wide tilde R from Equation (49) *)
R2[mu_, m_, n_] := (n - m)/(2 mu + 1) *
 ( NIntegrate[BernoulliB[2 mu + 1, x]/(m + x)^(2 mu + 1), {x, 0, 1/2},
     WorkingPrecision -> workPrec2] +
   NIntegrate[BernoulliB[2 mu + 1, x]/(m + x)^(2 mu + 1), {x, 1/2, 1},
     WorkingPrecision -> workPrec2] );


(* Equation (51) *)
lnFactRatioErrorBound[mu_, m_, n_] := Abs[R2[mu, m, n]] + Abs[R2[mu, n + 1, 2 n]] ;


lnFactRatio[mu_, m_, n_, iPrint_] :=
Module[
  (* compute the natural log of the ratio of two factorials:
       ln[ (2n)!/(n!)^2 ]; see Equation (50).
     also, throw away digits that are not justified based on R2[mu, m, n].
  *)
  { lnFactRat, lnFactErrorEst, lnFactRatShortened },
  (* note: we know initialLogSum to only (nDecimals2) decimal places.
      therefore, we cannot know lnFactRat to more decimal places than that.
      when Mathematica adds curlyPsi to initialLogSum, it will not retain
      more than (nDecimals2) decimal places in the result.
  *)
  lnFactRat = -initialLogSum - curlyPsi[mu, m, n] + curlyPsi[mu, n + 1, 2 n];
  lnFactErrorEst = lnFactRatioErrorBound[mu, m, n];
  lnFactRatShortened = removeQuestionableDigits[lnFactRat, lnFactErrorEst];
  If[iPrint == 1,
    Print["n = ", n, ", lnFactRat = ", lnFactRat, ", lnFactErrorEst = ", N[lnFactErrorEst, 5],
          ", lnFactRatShortened = ", lnFactRatShortened]
  ];

  Return[ lnFactRatShortened ]
]  (* end of Module *)




(* all necessary functions are now defined *)

iPrint = 0;  (* set to 1 to print intermediate data *)

(* expr1 and const1 are used only to display the expression used in FindRoot;
   in the actual calculations, we use curlyPhi[ ].
*)
If[iPrint == 1,
  expr1 = Simplify[ curlyPhi[mu1, m1, x] ];
  const1 = Part[expr1, 1];  (* extract the rational number *)
  Print["approximate phi to be used in FindRoot: ", (expr1 - const1) + N[const1, nDecimals1]];
];

(* calculate everything we need for this range of a0 values *)
a0First = 10;
a0Last = 25;

Print["integral for a0 is (Pi/2) * (1 - t)"];

a0List = nList = sumAkList = r1List = mantList = expoList = {};
For[a0 = a0First, a0 <= a0Last, a0++,
  { n , sumAk , r1 } = getNValueAndSumForA0[a0, m1, mu1, initialAkSum, accGoal1, workPrec1, iPrint];
  If[n == 0, Break[] ];  (* could not compute a valid value of n or sum of ak *)
  (* compute the natural log of (2n)!/(n!)^2 *)
  lnFactorialRatio = lnFactRatio[mu2, m2, n, iPrint];
  (* compute log of t to base 10; see Equations (38) - (40) *)
  logTBase10 = (n * Log[sumAk - a0] - (2 n - 1) * Log[2] + lnFactorialRatio) / Log[10];
  { mant , expo } = getME[ logTBase10 ];
  (* save lists of these values in case we want to use them later *)
  AppendTo[a0List, a0];
  AppendTo[nList, n];
  AppendTo[sumAkList, sumAk];
  AppendTo[r1List, r1];
  AppendTo[mantList, mant];
  AppendTo[expoList, expo];
  Print["a0 = ", a0, ", n = ", n, ", t = ", mant, " * 10^", expo];
]  (* end For a0 loop *)

\end{lstlisting}
\normalsize

\subsection{Results} \label{S:MMACode-Results}

Here are the values of the integrals for $a_0 = 10$ through $a_0 = 25$, obtained by running the code in Section \ref{S:MMACode-EulerNac}.
The value of the integral for $a_0$ is $(\pi/2) \cdot (1 - t)$.

If a computed intermediate value had more digits after the decimal point than were justified by the corresponding error value (for example, Equation \eqref{E:error_bound} or \eqref{E:factorialRatioErrorBound}), then the extra digits were discarded.
Therefore, all digits shown below should be correct, rounded in the last decimal place.

\tiny
\begin{lstlisting}
integral for a0 = (Pi/2) * (1 - t)
a0 = 10, n = 68100151, t = 9.649273600428684463479553120939810530923242208735398 * 10^-554381308
a0 = 11, n = 503195829, t = 7.57929806494947536128349934756162195412431861759227 * 10^-4887781043
a0 = 12, n = 3718142208, t = 5.30200436015724605246826614752108917188558325098544 * 10^-39227165565
a0 = 13, n = 27473561358, t = 1.52739916984845667363367296109645541442392755493153 * 10^-297230209953
a0 = 14, n = 203003686106, t = 6.617345077783595182168242992545965700461478406777 * 10^-2419966945909
a0 = 15, n = 1500005624924, t = 5.26019597269976433379615815051550875066124252042 * 10^-18988869014266
a0 = 16, n = 11083625711271, t = 4.06751521421327233190115950829638686972451899823 * 10^-148452517153987
a0 = 17, n = 81897532160125, t = 2.8703074957720537216132995767534053015103162770 * 10^-1261337931785960
a0 = 18, n = 605145459495141, t = 1.46932966274512803735093876340436661798499994 * 10^-9192758406970262
a0 = 19, n = 4471453748222757, t = 7.36887339695623805028019042180415757921528157 * 10^-73134639260589997
a0 = 20, n = 33039822589391676, t = 5.8024461422390775663611817270349845938468954 * 10^-579426465025122292
a0 = 21, n = 244133102611731231, t = 1.3869021709986676325063938918439160007102035 * 10^-4427143349945912840
a0 = 22, n = 1803913190804074904, t = 4.10151245193385022941804060193305405447094 * 10^-34064698104956009918
a0 = 23, n = 13329215764452299411, t = 1.71715972357092319138607947022428968556534 * 10^-259489336406929338805
a0 = 24, n = 98490323038288832267, t = 1.4982758030623762036996263232893870122517 * 10^-2011250066953860707590
a0 = 25, n = 727750522131718025058, t = 2.723848647528233561685563127899349771647 * 10^-15968197862152240928105
\end{lstlisting}
\normalsize

The above code began by initializing these parameters:
\begin{lstlisting}
  m1 = 100001, mu1 = 10, nDecimals1 = 100, accGoal1 = 20, workPrec1 = 40
  m2 = 100001, mu2 = 5, nDecimals2 = 100, workPrec2 = 40
\end{lstlisting}

The reader can obtain additional significant digits in $t$ by increasing some of these parameters.
Much of the time is already spent computing $s(m1-1)$, so it is more efficient to increase \verb+mu1+ or \verb+mu2+.

These calculations can also be extended beyond $a_0 = 25$.
For example, running the code with \verb+a0Last = 40+ with the above parameters gives, for $a_0 = 40$,
\SMALL
\begin{lstlisting}
t = 1.8758610 * 10^-266134053348172015148849587491648267
\end{lstlisting}
\normalsize
As $a_0$ increases, more and more of the significant digits in the calculation are consumed in the exponent.
Merely increasing \verb+mu1+ to 11 gives, for $a_0 = 40$,
\SMALL
\begin{lstlisting}
t = 1.8758609814211976 * 10^-266134053348172015148849587491648267
\end{lstlisting}
\normalsize

\end{document}